\documentclass[12pt]{amsart}
\usepackage{latexcad}
\usepackage{natbib}
\theoremstyle{definition}
\newtheorem{definition}{Definition}[section]

\theoremstyle{plain}
\newtheorem{theorem}[definition]{Theorem}
\newtheorem{lemma}[definition]{Lemma}

\theoremstyle{remark}

\newtheorem{remark}[definition]{Remark}

\def\to{\longrightarrow}

\def\span#1{\langle#1\rangle}

\def\neutral{e}

\def\chein{M_{2n}(G,\,2)}

\title{The Smallest Moufang Loop Revisited}
\author{Petr Vojt\v echovsk\'y}
\address{\scriptsize Department of Mathematics, University of Denver, 2360 S Gaylord St, Denver, CO  80208, USA}
\email{\scriptsize petr@math.du.edu
{\normalsize\begin{flushright}\textrm{Eingegangen am 13.\ M\"{a}rz
2003}\end{flushright}}}

\thanks{Work partially supported by Grant Agency of Charles University, grant
number 269/2001/B-MAT/MFF}

\begin{document}

\begin{abstract} We derive presentations for Moufang loops of type $\chein$,
defined by Chein, with $G$ finite, two-generated. We then use $G=S_3$ to
visualize the smallest non-associative Moufang loop.

\noindent {\sc Math.\ Subject Classification:} Primary: 20N05, Secondary: 20F05

\end{abstract}

\maketitle

\thispagestyle{empty}

\section{Introduction}

\noindent In order to derive a presentation for a groupoid $A=(A,\,\cdot)$, one
usually needs to introduce a normal form for elements of $A$ written in terms
of some generators. Such a normal form is not easy to find when $A$ is not
commutative, and even more so when $A$ is not associative. Once a normal form
is found, it might be still  difficult to come up with presenting relations.
Indeed, it is often the case that the only known presentation for a
non-associative groupoid is the \emph{table presentation}, i.e., the
presentation consisting of all relations $x\cdot y=z$ such that $x\cdot y$
equals $z$ in $A$, and where $x$, $y$ run over all elements of $A$. Table
presentations are extremely useful when one constructs a multiplication table
for $A$, however, they are of little use when one needs to identify $A$ as a
subgroupoid of another groupoid. To do the latter, it is necessary, in
principle, to evaluate all products $x\cdot y$ with $x$, $y\in A$. It is
therefore desirable to have access to presentations with a few presenting
relations.

The infinite class of Moufang loops $\chein$, defined below, represents a
significant portion of non-associative Moufang loops of small order. We derive
compact presentations for $\chein$ for every finite, two-generated group $G$.

Thirty years ago, Chein and Pflugfelder \cite{CheinPflugfelder} proved that the
smallest non-associative Moufang loop is of order $12$ and is unique up to
isomorphism. It coincides with $M=M_{12}(S_3,\,2)$. Guided by our presentation
for $M$, we give a new, visual description of $M$ in the last section.

\section{The Loops $M_{2n}(G,\,2)$}

\noindent A loop $L=(L,\,\cdot)$ is \emph{Moufang} if it satisfies one of the
three equivalent \emph{Moufang identities}
\begin{equation}\label{Eq:MI}
    xy\cdot zx=x(yz\cdot x),\;\;
    x(y\cdot xz)=(xy\cdot x)z,\;\;
    x(y\cdot zy)=(xy\cdot z)y.
\end{equation}
In fact, it is not necessary to assume that $L$ possesses a neutral element. By
a result of Kunen \cite{Kunen}, every quasigroup satisfying one of the Moufang
identities is necessarily a (Moufang) loop. Every element $x$ of a Moufang loop
has a two-sided inverse $x^{-1}$. Also, Moufang loops are \emph{diassociative},
i.e, every two-generated subloop is a group. We will use these well-known
properties of Moufang loops without warning throughout the paper.

The following construction is due to O.~Chein \cite{CheinII}. Let $G$ be a
finite group of order $n$. Pick a new element $u$, and define
\begin{displaymath}
    \chein=\{gu^{\alpha};\; g\in G,\,\alpha=0,\,1\},
\end{displaymath}
where
\begin{equation}\label{Eq:Mult}
    gu^{\alpha}\cdot hu^{\beta}
    =(g^{(-1)^{\beta}}h^{(-1)^{\alpha+\beta}})^{(-1)^{\beta}}
        u^{\alpha+\beta}\;\;\;(g,\,h\in G,\,\alpha,\,\beta=0,\,1).
\end{equation}
Then $\chein$ is a Moufang loop of order $2n$. It is associative if and only
if $G$ is commutative.

Let $\pi(m)$ be the number of isomorphism types of non-associative Moufang
loops of order at most $m$, and let $\sigma(m)$ be the number of
non-associative loops of the form $\chein$ of order at most $m$. Then,
according to Chein's classification \cite{CheinII}, $\pi(31)=13$,
$\sigma(31)=8$, $\pi(63)=158$, $\sigma(63)=50$. (As Orin Chein kindly notified
me, Edgar Goodaire noticed that the loop $M_{12}(S_3,\,2)\times C_3$ is missing
in \cite{CheinII}. He also observed that $M_{48}(5,\,5,\,5,\,3,\,3,\,0)$ is
isomorphic to $M_{48}(5,\,5,\,5,\,3,\,6,\,0)$, and
$M_{48}(5,\,5,\,5,\,3,\,3,\,6)$ to $M_{48}(5,\,5,\,5,\,3,\,6,\,6)$. That is
why $\pi(63)$ equals $158$, rather than $159$.) This demonstrates eloquently
the abundance of loops of type $\chein$ among Moufang loops of small order.

\section{The Presentations}

\noindent We start with the table presentation $(\ref{Eq:Mult})$ for $\chein$
and prove

\begin{theorem}\label{Th:Main}
Let $G=\span{x,\,y;\;R}$ be a presentation for a finite group $G$, where $R$
is a set of relations in generators $x$, $y$. Then $\chein$ is presented by
\begin{equation}\label{Eq:Pres}
    \span{x,\,y,\,u;\; R,\,u^2=(xu)^2=(yu)^2=(xy\cdot u)^2=\neutral},
\end{equation}
where $\neutral$ is the neutral element of $G$.
\end{theorem}

Let us emphasize that $(\ref{Eq:Pres})$ is a presentation in the \emph{variety
of Moufang loops}, not groups.

The complicated multiplication formula $(\ref{Eq:Mult})$ merely describes the
four cases
\begin{eqnarray}
    g\cdot h &=& gh,\label{Eq:Case1}\\
    gu\cdot h &=& gh^{-1}\cdot u,\label{Eq:Case2}\\
    g\cdot hu &=& hg\cdot u,\label{Eq:Case3}\\
    gu\cdot hu &=& h^{-1}g\label{Eq:Case4}
\end{eqnarray}
in a compact way. In particular, identities $(\ref{Eq:Case4})$ and
$(\ref{Eq:Case2})$ imply
\begin{equation}\label{Eq:ShortPres}
    u^2=\neutral,\;\; gu=ug^{-1}\ (g\in G).
\end{equation}
We claim that $(\ref{Eq:ShortPres})$ is equivalent to $(\ref{Eq:Mult})$. An
element $g\in G$ will be called \emph{good} if $gu=ug^{-1}$ can be derived
from $(\ref{Eq:Pres})$.

\begin{lemma}\label{Lm:Equiv}
If $h\in G$ is good, then $(\ref{Eq:Case2})$ holds. If $g$, $h$, $hg\in G$ are
good, then $(\ref{Eq:Case3})$ holds. If $g$, $g^{-1}h$ are good, then
$(\ref{Eq:Case4})$ holds.
\begin{proof}
We have $gu\cdot h=(gu\cdot h)u\cdot u = (g\cdot uhu)u = (g\cdot h^{-1}uu)u =
gh^{-1}\cdot u$ if $h$ is good. Assume that $g$, $h$, $hg$ are good. Then
$g\cdot hu = g\cdot uh^{-1} = u\cdot u(g\cdot uh^{-1}) = u(ugu\cdot h^{-1}) =
u\cdot g^{-1}h^{-1}= hg\cdot u$. Finally, when $g$ and $g^{-1}h$ are good, we
derive $gu\cdot hu = ug^{-1} \cdot hu = u\cdot g^{-1}h\cdot u = h^{-1}g$.
\end{proof}
\end{lemma}

Thus $(\ref{Eq:ShortPres})$ is equivalent to $(\ref{Eq:Mult})$. Moreover, in
order to prove Theorem \ref{Th:Main}, it suffices to show that every $g\in G$
is good.

Thanks to diassociativity, $g^{s}$ ($s$ positive integer) is good whenever $g$
is. Since $G$ is finite, $g^{-1}$ is good whenever $g$ is.

\begin{lemma}\label{Lm:Swap}
Assume that $g$, $h\in G$ are good. Then $gh$ is good if and only if $hg$ is.
\begin{proof}
Because of the symmetry, it is enough to prove only one implication. Assume
that $hg$ is good. By Lemma \ref{Lm:Equiv}, $g\cdot hu = hg\cdot u$. Using this
identity, we obtain $gh\cdot ug = g(hu\cdot g) = (g\cdot hu)g = (hg\cdot u)g =
h\cdot gug = hu$, thus $gh = hu\cdot g^{-1}u = uh^{-1}\cdot g^{-1}u = u\cdot
h^{-1}g^{-1}\cdot u$, and so $gh\cdot u = u\cdot h^{-1}g^{-1}$.
\end{proof}
\end{lemma}

\begin{lemma}\label{Lm:Up}
Assume that $g$, $h\in G$ are good. Then so is $ghg$.
\begin{proof}
Since $g^{-1}$, $h$ are good, Lemma \ref{Lm:Equiv} yields $ug\cdot h =
g^{-1}u\cdot h = g^{-1}h^{-1}\cdot u$. Then $u\cdot ghg\cdot u = (ug\cdot
h)g\cdot u = (g^{-1}h^{-1}\cdot u)g\cdot u = g^{-1}h^{-1}\cdot ugu =
g^{-1}h^{-1}g^{-1}$, and we are done.
\end{proof}
\end{lemma}

We continue by induction on the \emph{complexity}, or \emph{length}, if you
will, of the elements of $G$, defined below.

For $\varepsilon=1$, $-1$, let $X_\varepsilon$ be the set of symbols
$\{x_1^\varepsilon$, $\cdots$, $x_m^\varepsilon\}$, and write $X=X_1\cup
X_{-1}$. Every word $w$ of the free group $F=\span{X}$ can be written uniquely
in the form $x_{i_1}^{\varepsilon_1}\cdots x_{i_r}^{\varepsilon_r}$, where
$i_j\ne i_{j+1}$, and $\varepsilon_j$ is a non-zero integer. Define the
\emph{complexity} of $w$ as the ordered pair
$c(w)=(r,\,\sum_{j=1}^r|\varepsilon_r|)$, and order the complexities
lexicographically.

From now on, assume that $G$ is two-generated, and write $x=x_1$, $y=x_2$.

Since $xu=ux^{-1}$ and $yu=uy^{-1}$ are presenting relations, both $x$, $y$ are
good, and hence both $x^s$, $y^s$ are good for every integer $s$. The last
presenting relation $xy\cdot u=u\cdot y^{-1}x^{-1}$ shows that both $xy$ and
$y^{-1}x^{-1}=(xy)^{-1}$ are good. Then $yx$ and $x^{-1}y^{-1}=(yx)^{-1}$ are
good, by Lemma \ref{Lm:Swap}. Also, Lemma \ref{Lm:Up} implies that
$x^{-1}\cdot xy\cdot x^{-1}=yx^{-1}$ is good. Then $x^{-1}y$,
$xy^{-1}=(yx^{-1})^{-1}$ and $y^{-1}x = (x^{-1}y)^{-1}$ are good, by Lemma
\ref{Lm:Swap}. This means that every $g\in G$ with $c(g)<(2,\,3)$ is good.

\begin{lemma}\label{Lm:MiddleStep}
Every $g\in G$ with $c(g)<(3,\,0)$ is good.
\begin{proof}
Suppose there is $g$ that is not good, and let $c(g)=(r,\,s)$ be as small as
possible. We can assume that $g=a^ub^v$, where $\{a,\,b\}=\{x,\,y\}$,
$s=|u|+|v|>2$, and $u\ne 0\ne v$.

Either $|u|>1$ or $|v|>1$. Without loss of generality, $u>1$. (By Lemma
\ref{Lm:Swap}, we can assume that $|u|>1$. When $u$ is negative, consider the
inverse $b^{-v}a^{-u}$ instead, and apply Lemma \ref{Lm:Swap} again.) Since
$c(a^{u-2}b^v)<(2,\,s)$, the element $a^{u-2}b^v$ is good, and so is
$a^{u-1}b^va = a\cdot a^{u-2}b^v\cdot a$. As $a^{u-1}b^v$ is good by the
induction hypothesis, $a^ub^va = a\cdot a^{u-1}b^v\cdot a$ is good as well, by
Lemma \ref{Lm:Up}. Then the decomposition of $a^{u-1}b^va$ into $a^{-1}\cdot
a^ub^va$ demonstrates that $a^ub^va\cdot a^{-1}=a^ub^v$ is good, by Lemma
\ref{Lm:Swap}. We have reached a contradiction.
\end{proof}
\end{lemma}

To finish the proof, assume there is $g\in G$ that is not good, and let
$c(g)=(r,\,s)$ be as small as possible. By Lemma \ref{Lm:MiddleStep}, $r\ge 3$.
When $r$ is odd, we can write
$g=a^{\varepsilon_1}b^{\varepsilon_2}a^{\varepsilon_3}\cdots
b^{\varepsilon_{r-1}}a^{\varepsilon_r} = khk$, where $k=a^{\varepsilon_r}$,
$h=a^{\varepsilon_1-\varepsilon_r}b^{\varepsilon_2}a^{\varepsilon_3} \cdots
b^{\varepsilon_{r-1}}$, and $\{a,\,b\}=\{x,\,y\}$. Since $c(k)$,
$c(h)<(r,\,s)$, both $k$, $h$ are good, and then $g$ is good by Lemma
\ref{Lm:Up}.

Assume that $r$ is even. Then $g=a^{\varepsilon_1}b^{\varepsilon_2}\cdots
a^{\varepsilon_{r-1}}b^{\varepsilon_r} = khk$, where $k=a^{\varepsilon_1}
b^{\varepsilon_r}$ and $h=b^{\varepsilon_2-\varepsilon_r}a^{\varepsilon_3}
\cdots b^{\varepsilon_{r-2}}a^{\varepsilon_{r-1}-\varepsilon_1}$. Again,
$c(k)$, $c(h)<(r,\,s)$, thus both $k$ and $h$ are good, and so is $g$, by Lemma
\ref{Lm:Up}.

Theorem \ref{Th:Main} is proved.

\section{Visualization of the Smallest Moufang Loop}

\noindent The multiplication formula $(\ref{Eq:Mult})$ for $M=M_{12}(S_3,\,2)$
is certainly difficult to memorize, and so is the one in \cite[Example
IV.1.2]{Pflugfelder}. We present a visual description of $M$.

Note that there are $9$ involutions and $2$ elements of order $3$ in $M$ (cf.
\cite[Table 3]{CheinI}). We are going to define a $12$-element groupoid $L$ and
show that it is isomorphic to $M$.

\begin{figure}[ht]
    \centering
    \input{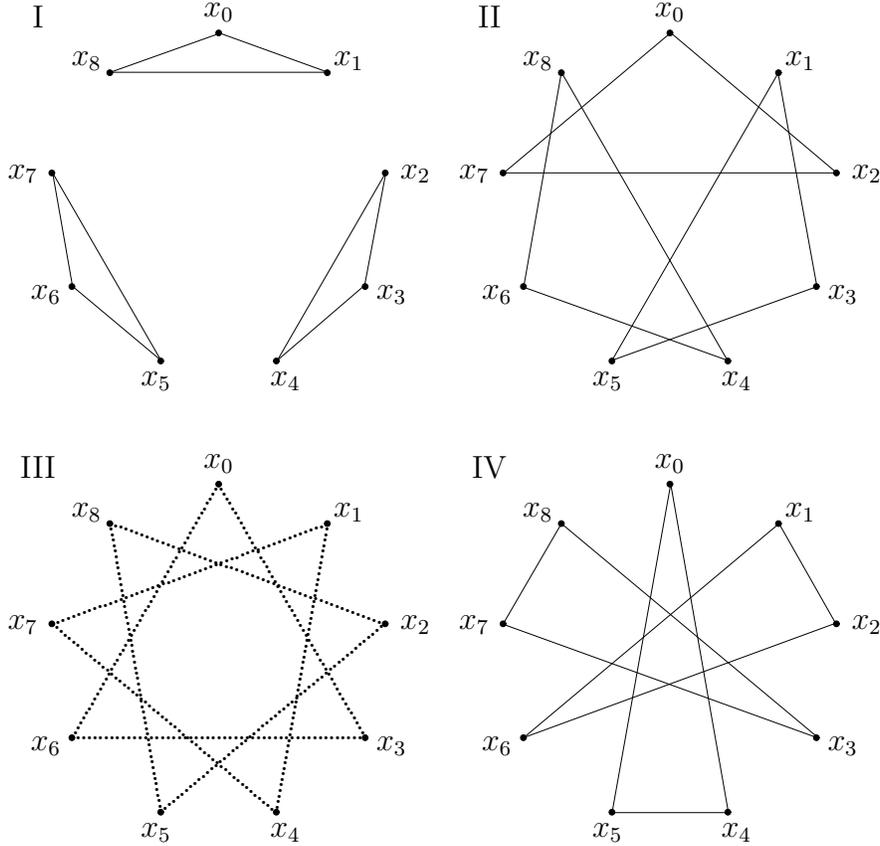}
    \caption[]{Multiplication in $M_{12}(S_3,\,2)$}
    \label{Fg:M12}
\end{figure}

Look at the four diagrams in Figure \ref{Fg:M12}. Think of the vertices $x_0$,
$\dots$, $x_8$ as involutions. Let $L$ consists of $e$, $x_0$, $\dots$, $x_8$,
$y$, $y^{-1}$, where $y$ is of order $3$. Interpret the edges of diagrams I--IV
as multiplication rules in the following way. If $x_i$ and $x_j$ are connected
by a solid line, let $x_ix_j$ be the third vertex of the (unique) triangle
containing both $x_i$ and $x_j$. If $x_i$ and $x_j$ are not connected by a
solid line, we must have $j=i\pm 3$, and then $x_i$ and $x_j$ are connected by
a dotted line (in diagram III). Define $x_ix_{i+3}=y$ and $x_ix_{i-3}=y^{-1}$.

\setlength{\unitlength}{0.75mm}

This partial multiplication can be extended by properties of Moufang loops. To
avoid ambiguity, we postulate that $x_iy=y^{-1}x_i=x_{i+3}$ and
$yx_i=x_iy^{-1}=x_{i-3}$.

Obviously, $L$ is closed under multiplication and has a neutral element. It is
non-associative, since $x_0x_1\cdot x_3=x_8x_3=x_7\ne x_4=x_0x_5=x_0\cdot
x_1x_3$. Is $L$ isomorphic to $M$? There is a unique Moufang loop of order $12$
\cite{CheinPflugfelder}, so it suffices to check the Moufang identities for
$L$. However, this is not so easy! Instead, we verify directly that $L$
satisfies the multiplication formula $(\ref{Eq:Mult})$ with some choice of $G$
and $u$.

\begin{remark}
It does not suffice to verify $(\ref{Eq:ShortPres})$ for some choice of $G$
and $u$ because $(\ref{Eq:ShortPres})$ is equivalent to $(\ref{Eq:Mult})$ only
when it is assumed that $L$ is Moufang.
\end{remark}

Put $x=x_0$, and observe that $G=\span{x,\,y}=\{\neutral$, $x_0$, $y$, $x_3$,
$x_6$, $y^{-1}\}$ is isomorphic to $S_3$. Let $u=x_1\not\in G$. We show that
$(\ref{Eq:Case1})$--$(\ref{Eq:Case4})$ are satisfied for every $g$, $h\in G$.
Thanks to the symmetry of Figure \ref{Fg:M12}, it is enough to consider only
$\{g,\,h\}=\{x_0,\,x_3\}$, $\{x_0,\,y\}$.

Identity $(\ref{Eq:Case1})$ is trivial. Let us prove $(\ref{Eq:Case2})$. We
have $x_0x_1\cdot x_3=x_8x_3 = x_7 = yx_1 = x_0x_3^{-1}\cdot x_1$, $x_0x_1\cdot
y = x_8y = x_2 = x_6x_1 = x_0y^{-1}\cdot x_1$, $x_3x_1\cdot x_0 = x_5x_0 = x_4
= y^{-1}x_1 = x_3 x_0^{-1}\cdot x_1$, and $yx_1\cdot x_0 = x_7x_0 = x_2 =
x_6x_1 = yx_0^{-1}\cdot x_1$. Similarly for $(\ref{Eq:Case3})$,
$(\ref{Eq:Case4})$.

Hence $L$ is isomorphic to $M$. The subloop structure of $L$ is apparent from
the visual rules, too. If $j\equiv i\pmod 3$ then $\span{x_i,\,x_j}\cong S_3$;
otherwise, $\span{x_i,\,x_j}\cong V_4$, for $i\ne j$.

\bibliographystyle{plain}

\end{document}